\documentclass[12pt]
{article}
\usepackage{amssymb,amsmath, mathtools}

\usepackage[latin1]{inputenc}
\usepackage[T1]{fontenc}
\usepackage{play}

\usepackage{graphicx}
\DeclareGraphicsExtensions{.pdf,.png,.jpg}

\graphicspath{ {images/} }

\usepackage{lipsum}

\let\OLDthebibliography\thebibliography
\renewcommand\thebibliography[1]{
  \OLDthebibliography{#1}
  \setlength{\parskip}{0pt}
  \setlength{\itemsep}{0pt plus 0.3ex}
}

\usepackage{amsmath}
\usepackage{amsthm}
\usepackage{amsfonts}
\usepackage{bbm}
\usepackage{amssymb}
\usepackage{graphicx}
\usepackage{color}

\usepackage{wasysym, marvosym}

\usepackage{fancybox}
\newlength{\querylen}
\newtheorem{thm}{Theorem}
\newtheorem{lemma}[thm]{Lemma}

\newtheorem{cor}[thm]{Corollary}
\newtheorem{example}[thm]{Example}
\newtheorem{assertion}[thm]{Proposition}
\theoremstyle{definition}
\newtheorem{define}[thm]{Definition}
\theoremstyle{remark}

\begin{document}
\title{How to beat the $1/e$-strategy of best choice\\
(the random arrivals problem)} 
\author{Alexander Gnedin\thanks{Queen Mary, University of London}}
\maketitle

\begin{abstract}
\noindent
In the best choice problem with random arrivals,
an unknown number  $n$ of rankable items arrive at times sampled from the uniform distribution.
As is well known,  a real-time player can ensure stopping at the  overall best item with probability at least $1/e$, by waiting until 
time $1/e$ then selecting the first relatively best item to appear (if any).
This paper discusses the issue of  dominance in a wide class of stopping strategies of best choice, and argues that 
in fact the player  faces a trade-off between success probabilities for various values of $n$.
We show that the $1/e$-strategy is not a unique minimax strategy  and that it can be improved in various  ways.
\end{abstract}

\section{Introduction} 
In a version of the familiar
 best choice problem,  $n$ items ranked $1,2,\cdots,n$ from best to worst
 arrive at independent,  uniformly distributed random times on $[0,1]$. 
Alice, a real-time
player, observes  arrivals in the chronological order, and 
 ranks each  item
 relative to all other seen so far.
She can stop  anytime, winning if  the stopping occurs at  the best of $n$ items and losing in any other case.
The question is about good stopping strategies when
the player does not know $n$.

An item  is said to be a record if  there is no better item that arrived before it. 
The first observed item is a record,  and also any subsequent arrival with relative rank one. 
The overall best item appears as the last record, therefore 
strategies ever stopping on non-records should be discarded straight away.

If  Alice knew $n$ she could  achieve
the maximum  winning probability using  the classic {\it $d^*$-strategy}:
wait until $d^*-1$  items pass by  (where $d^*=d^*(n)$ is about $n/e$)  then
stop at the first subsequent 
record, if available. 
For every $n$, 
the winning   probability of the {$d^*$-strategy} exceeds $1/e$, and approaches
the lower bound  monotonically as $n$ increases. By certainty about $n$ 
observing the arrival times is of no avail, because by exchangeability they  
are independent of the relative ranks.


However, if $n$  is unknown,  
no strategy based only on  counting and relative ranking of the arrivals can ensure a winning  probability bounded away from zero for all $n$ \cite{Abdel}. 
Starting from \cite{Presman} such strategies were studied under the assumption that
$n$ is drawn from some distribution known or partly known to the player, 
or under  an upper constraint on the number of items  \cite{HillKrengel, HillKennedy}. 

At first glance if the player is completely  ignorant about $n$  there is no useful alternative, but 
 this is where the time factor steps in.  For $n$ large, the number of items arrived by time $t$ is close to $nt$, thus the {\it $1/e$-strategy}, which prescribes
to wait until time $1/e$ then stop at the first subsequent record (if any), will have about the same  effect as the $d^*$-strategy that uses the full knowledge of $n$. 
Intuitively, 
in the random arrivals model, 
the  lack of information about $n$ is compensated 
by the proportionate growth of the number of observations.

Looking back in history,  the $1/e$-strategy seems to have made its first appearance
 in  \cite{Gaver}, where asymptotic optimality was noticed for the model with  items observed at times of the  Poisson process.
The  exact  optimality was shown in Rubin's `secretary problem' with infinitely many arrivals \cite{BG, GianiniSam, Rubin}, 
and in a Bayesian  setting with $n$  sampled from the improper uniform distribution on integers \cite{Stewart},
but  such  limit forms of the problem are very different in spirit and will not be touched upon here.
Stewart \cite{Stewart}, assuming i.i.d. exponentially distributed arrival times, proved asymptotic optimality  of the strategy with  the
cutoff time chosen to be the $1/e$ quantile of the exponential distribution, and showed numerical evidence that the strategy performs well  also for small values of $n$.
Bruss \cite{Bruss84} gave the strategy its name and made the key observations that 
 the distribution of arrival times is not important
(provided it is continuous and known to the player) and
that the benchmark winning  probability $1/e$ appears as the sharp lower bound,  approached monotonically as $n$ increases.

The findings from \cite{Bruss84} 
taken together with the upper bound  of  the $d^*$-strategy imply that the $1/e$-strategy is minimax. This naturally  begs two questions:
is the $1/e$-strategy unique minimax, 
and is there something better? 
Corollary 2 from \cite{Bruss84}   stated the  uniqueness within a class of strategies that wait certain fixed time then stop at the $r$th record to follow
(though second part of the definition of the class on p. 883 is controversial).
Subsequent work  citing the result upgraded it to a proof of the overall uniqueness of minimax strategy,   
and 
added further optimality properties, although no rigorous arguments were given (see \cite{Bruss88} p. 313, \cite{BrussYor} p. 3259).
This made the way to surveys and popular  sites forming a consensus in the mathematical community that 
no strategy could be better. 
In this paper we disprove these uniqueness and optimality claims from the position of game-theoretic dominance. Concretely, we show that
\begin{itemize}
\item[(i)]  the $1/e$-strategy is never optimal if $n$ is drawn from some probability distribution known to Alice,
\item[(ii)] there are many  strategies that achieve the bencmark $1/e$ asymptotically,
\item[(iii)] for every  $u$ there exist strategies that strictly improve  the $1/e$-strategy 
simultaneously for all $n\leq u$,
\item[(iv)]  there exists a simple strategy outperforming the $1/e$ strategy simultaneously for all $n>1$
(strictly for $n>2$),
\item[(v)] there exist more complex strategies strictly outperforming the $1/e$ strategy simultaneously for all $n>2$,
\item[(vi)] for every $\ell\geq 1$ there exist still more complex  strategies that guarantee 
the winning probability at least $1/e$ for all $n$, and are outperforming
the $1/e$-strategy simultaneously for all $n>\ell$.
\end{itemize}
In practical terms, for $\ell=1$  Alice faces a trade-off between some
advantage when the  choice  occurs from  many opportunities against a higher risk of
 going away empty-handed  if  just a sole item shows up.

A reason for stating (iv)  with general  $\ell\geq 1$, rather than   just $\ell=1$, is that the dominance can be  also explored 
in the opposite direction. In the limit $\ell\to\infty$  we shall arrive at the following finding:
\begin{itemize}
\item[(vii)]  there exist minimax strategies (hence winning with chance above $1/e$), which are worse than the $1/e$-strategy simultaneously for  all $n$.
\end{itemize}
Although the  game  has a value (see Corollary 3 below), by (i) there is no worst-case distribution for the number of items, 
which paved the way to (vi) and (vii).
This extends the list of known dominance paradoxes  in the optimal stopping games of best choice \cite{googol, Guess, GnedinKr}. 
A question which we find difficult and leave here open is the existence of  unconstrained dominating strategy  ($\ell=0$ in (vi)).
If  the $1/e$-strategy turns to be undominated,  (vi) would mean that  the source of rigidity lies in the trivial case $n=1$.

We  shall make use of strategies that  rely  decisions at record times also on the number of arrivals seen so far.
Such strategies, called in \cite{BrussSam2} `nonstationary',
have been used in selection problems with random number of items \cite{ Browne, CZ, Ano, Tamaki}.
Exploring this large class, we  further construct  counter-examples to the assertion of 
 risk monotonicity  relative to the stochastic ordering on distributions of the number of items (see Theorem 2.3 in \cite{BrussSam1}, Equation (36) in \cite{Bounds},
Equation (76) in \cite{Handbook}).
Failing monotonicity of the kind was observed long ago in the rank-based  optimal stopping  problems with random sample size and fixed arrival times \cite{Gianini-P}.

\section{The game strategies}

Following a suggestion from \cite{Bruss84}, 
we adopt here  the paradigm of
 game theory, hence assuming that an antagonist player, Pierre, is  in charge of the variable $n$.

Pierre's strategies are easy to describe. A pure strategy is  just a positive integer, and a mixed strategy is
  a  probability  distribution $\nu=(\nu_1,\nu_2,\cdots)$.

When Pierre plays $\nu$, Alice sees the items coming at epochs of a random counting process $N=(N_t, ~t\in[0,1])$, where $N_1$ follows distribution $\nu$.
By the assumption of uniform arrival times, $N$ 
 has the order statistics property: conditional on $N_{t-}=k$  the configuration of  the first $k$ ordered arrivals  is uniformly distributed on the $k$-simplex
$\{(t_1,\cdots, t_{k}): 0<t_1<\cdots<t_{k}<t\}$.
Equivalently,  $N$ can be regarded as  a Markovian, inhomogeneous birth process with $N_0=0$ and the jump rates readily 
computable in terms of the probability generating function of $\nu$, see  \cite{Puri}.
The observed $k$th arrival  is ranked $R_k$th relative to the items seen so far, where the relative ranks $R_1, R_2,\cdots$ are independent and independent of $N$, with $R_k$ being uniform on $\{1,\cdots,k\}$. 
Thus the  process counting  record times is  derived from  $N$ by thinning, such that the $k$th arrival is kept 
with probability $1/k$ independently of anything else \cite{BrowneBunge}.

\paragraph{Remark} The best choice problem with arrivals coming by a point process is a 
well studied topic. To adjust a bulk of  past work \cite{Browne, BrussRogers, CZ, Gaver,  Tamaki} to the setting of this paper,
the processes need to be conditioned on  nonzero number of arrivals.
Equivalently, we may allow for
 extended
 distributions 
  on nonnegative integers, of the form $(1-c, \,c\nu_1, \,c\nu_2,\cdots)$, $0<c\leq1$.
Whatever $c$,  from  the first arrival on the conditional distribution of $N_1$ is the same.
For instance, in the setting of arrivals by Poisson process, the game against the Poisson distribution 
is  the same as against  the zero-truncated Poisson distribution.

\vskip0.2cm

The space of  Alice's  strategies is immense. 
The generic pure strategy $\tau$ is a stopping time, which is   as a certain choice (if any) from the set of record times.
 The eventually observed sample data includes  $n$ and a sequence $(t_1,r_1),\cdots, (t_n,r_n)$, where $t_k\in (0,1)$ are the increasing arrival times and     $r_k\in \{1,\cdots,k\}$  the relative ranks of items. The function $\tau$ chooses one of the $t_k$'s or $1$ (no stop) subject to the following conditions.
Firstly, if $t_k$ is chosen then the same choice is valid for any other $n'\geq k$ and possible `future' data
$(t_{k+1}', r_{k+1}'),\cdots,(t_{n'+1}', r_{n'+1}')$.
Secondly, $t_k$ can only be chosen if $r_k=1$ (record arrival).

A win with $\tau$ occurs in the event that $\tau=t_k$ (so, $r_k=1$) and $r_k>1,\cdots,r_{k+1}>1$ for some $k$.
Let $W(\tau,\nu)$ denote the winning probability when Alice plays stopping strategy $\tau$, and write $W(\tau,n)$ if Pierre's strategy is pure,
thus 
$$W(\tau,\nu)=\sum_{n\geq 1} \nu_n W(\tau,n).$$

\begin{define}
Stopping strategy  $\tau$ is called 

\begin{itemize}
\item[(a)] 
$n$-optimal if it achieves $\max W(\,\cdot\,,n)$,
\item[(b)]  $\nu$-optimal if  it achieves $\max W\,(\cdot\,,\nu)$,
\item[(c)] asymptotically optimal if $W(\tau, n)\to 1/e$ as $n\to\infty$,
\item[(d)] dominated  (respectively, dominated on a subset of positive integers $Z$)  if  there is  
another strategy $\tau'$ with $W(\tau',n)\geq W(\tau,n)$  for  all $n$ (respectively, for $n\in Z$), where at least one of the inequalities is strict,
\item[(e)] strongly  dominated on $Z$ if $W(\tau',n)>W(\tau,n), ~n\in Z$.
\end{itemize}
\end{define}
\noindent
In statistical decision theory, undominated strategy is also called  `admissible' \cite{Blackwell}.

Dicrete-time theory of optimal stopping ensures existence of maximisers in (b):
these are best-response counter-strategies of Alice.
By the nature of  arrivals process $N$, 
the sufficiency principle allows one to seek for maximisers within a reduced class of {\it Markovian} stopping times,
which make  decision at record  time $t$  only with the account of $t$ and the number of arrivals $N_t$. 
For $t$ is a record time, we call the number of observed items $N_t$  the {\it index} of  the record.
In particular, the earliest arrival is a record of index $1$.

We will restrict further consideration to Markovian strategies of special form
$$\tau=\min\{t:~ N_{t}-N_{t-}=1, ~R_{N_t}=1,~  N_t\geq a_{N_t}\}$$
($\min\varnothing=1$), where  $0\leq a_k\leq 1, ~k=1,2,\cdots$. The {\it cutoff}  $a_k$ specifies the 
 earliest time when a record with index $k$ can be accepted.

\subsection{$d$-strategies}

The classic {\it $d$-strategy}   stops at the first record of index at least $d$
regardless of the arrival time.
This has only extreme cutoffs
$$a_k=\begin{cases} 1,~~~k<d,\\
0,~~~k\geq d. 
\end{cases}
$$
The winning probability on $n$ items is  
\begin{eqnarray}\label{hh}    
s(d,n) =\frac{d-1}{n}\,h(d,n) ~{\rm with}~   ~h(d,n):= \sum_{j=d}^n\frac{1}{j-1},~~~d>1,
\end{eqnarray}
and $s(1,n)=1/n$. 
The $n$-optimal strategy has
\begin{equation}\label{d*}
d^*=\min\{d: h(d, n)\geq 1\},
\end{equation}
as is well known. For $n\neq2$ the maximiser 
is unique. 
For $n=2$, $s(1,2)=s(2,2)=1/2$; 
the non-uniqueness in this case leads to an interesting dominance phenomenon  related
to the Blackwell-Hill-Cover  paradox of guessing the larger of two numbers with probability exceeding $1/2$
\cite{googol, Guess}.

\subsection{$x$-strategies}

For $0\leq x\leq 1$ the strategy with identical cutoffs $a_k\equiv x$,
$$\tau_x:= \min\{t: ~t\geq x, ~R_{N_t}=1\},$$
stops at the first record arriving after time $x$ regardless of the index.
In particular,  $\tau_{1/e}$ is the $1/e$-strategy.

The winning probability is a polynomial $p_n(x):=W(\tau_x, n)$, which has 
several useful representations
\begin{eqnarray*}
p_n(x)=
\sum_{k=1}^{n-1} {n\choose k} x^k(1-x)^{n-k}s(k+1,n)&=&\\
\sum_{k=1}^{n-1} \frac{x(1-x)^k}{k}+\frac{(1-x)^n}{n}&=&\\
1-x-\sum_{k=2}^{n} \frac{(1-x)^{k}}{k(k-1)}.
\end{eqnarray*}
The first is obtained by conditioning on $N_x=k$ and  is analogous to Stewart's formula  for exponential arrivals \cite{Stewart}.
The second is obtained by conditioning on the event that $k$ top items arrive after $x$
and  the $(k+1)$st before $x$
\cite{Bruss84}.

The third can be argued by coupling over the sample size, as follows.
Suppose absolute ranks $1,\cdots, n-1$ appear as a winning (respectively, losing) configuration in the problem with $n-1$ items. 
Adding the worst $n$th arrival does not change the status quo unless
none of the arrivals falls in $[0,x]$, the best is the first appearing after $x$, and the worst falls between 
$x$ and the best arrival.  The formula follows then from the recursion
$$p_n(x)=p_{n-1}-\frac{(1-x)^n}{n(n-1)}.$$.

\begin{figure}\label{Figure1}
\includegraphics[width=9cm, height=7cm]{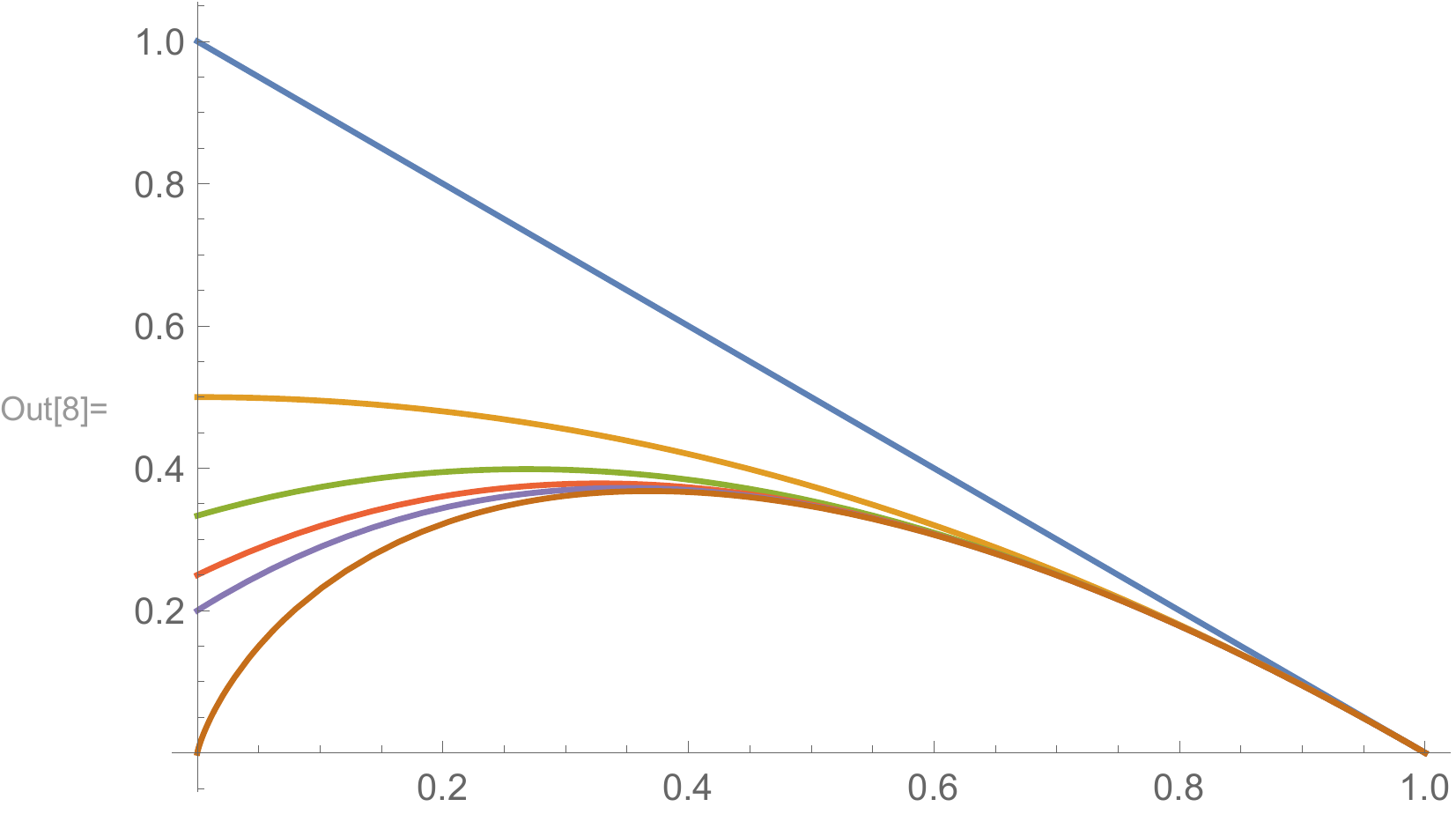}
\caption{$p_n(x)$ for $n=1,2,3,4,5$ and  $-x\log x$} 
\end{figure}

Note that $p_n(0)=1/n$.
From the third representation it is most easy to see that \cite{Bruss84}

\begin{itemize}
\item[(a)] $p_n(x)\downarrow -x\log x$ as $n\uparrow$,
\item[(b)] the $p_ns$'s are unimodal on $[0,1]$, with maximum points  $x_n$  satisfying  
$$0=x_1=x_2<x_3<\cdots \to 1/e,$$
\item[(c)] $p_1$ and $p_2$ are strictly decreasing.
\end{itemize}

\subsection{The case of geometric distribution}

If  Pierre plays  the geometric distribution $\nu_n=\theta (\theta+1)^{-n}$ (so with `success' probability $\theta/(\theta+1)$, ~$\theta>0$), the $\nu$-optimal strategy is $\tau_x$ with 
\begin{equation}\label{alpha-x}
x=\left(\frac{1}{e} -\frac{\theta (e-1)}{e}\right)_+,
\end{equation}
see \cite{Browne, Bruss-CZ, BrussSam2}. 
For $\theta\geq 1/(e-1)$ the winning probability with (\ref{alpha-x}) is $W(\tau_x,\nu)=(\theta+1)/e$. 
Note that $x=0$ for $\theta\geq 1/(e-1)$.

In this game,  Alice observes the items arriving at times of a random process, which is 
 mixed Poisson  with  random rate sampled from the exponential-$\theta$ distribution, and conditioned on nonzero number of arrivals.
The earliest arrival time $T$ has density $t\mapsto \theta(\theta+1)/(\theta+t)^2$,
and from this time on $N$ behaves like 
 a  pure-birth  process with birth rate 
$${\mathbb P}(N_{t+{\rm d}t}-N_t=1\,|N_t=k, t\geq T)=  \frac{k+1}{\theta+t}.$$
Since the $(k+1)$st arrival is a record with probability $1/(k+1)$, the cancellation of the state variable shows that the record-counting process on $[T,1]$ has deterministic compensator hence
is a Poisson process by Watanabe's theorem. 
This  fact underlies optimality of the $x$-strategy that disregards the index of record.
See  \cite{Browne,  BrussRogers, BrussSam2} for characterisation of this `stationary' case.

\section{Dominance and optimality}

Basic dominance features of the $x$-strategies 
 are drawn straight from the monotonicity properties of $p_n$'s.

\begin{assertion} The $x$-strategies satisfy

\begin{itemize}
\item[\rm(i)]  $\tau_x$ dominates $\tau_y$ for $1/e\leq x<y\leq 1$,
\item[\rm(ii)]  $\tau_x$ is dominated by $\tau_y$ on $\{n, n+1,\cdots\}$  for $0\leq x<y\leq x_n$,
\item[\rm (iii)] $\tau_x$ is undominated for $0\leq x<1/e$.
\end{itemize}
\end{assertion}

\begin{proof} The first two assertions are clear from the properties of $p_n$'s.
The third follows from $\nu$-optimality, for $\nu$ the geometric distribution with parameter found from (\ref{alpha-x}).
\end{proof}

\begin{cor} The stopping game has the following features:
\begin{itemize}
\item[\rm(i)] the value of the game is $1/e$, 
\item[\rm(ii)] Pierre has no minimax strategy and so the game has no saddle-point, 
\item[\rm(iii)] the $1/e$-strategy is never $\nu$-optimal.
\end{itemize}
\end{cor}
\begin{proof} We have $W(\tau_{1/e},\nu)>1/e$, but Pierre can keep the winning probability below $1/e+\varepsilon$ by playing
large $n$ or a stochastically large $\nu$. This gives (i).

Since $p_n(1/e)>1/e$, also $W(\tau_{1/e},\nu)>1/e$ for any $\nu$, whence (ii).

To argue (iii), note that
$p_n'(1/e)<0$ implies $\sum_n \nu_n p_n'(1/e)<0$, thus $x\mapsto W(\tau_x,\nu)$ is strictly decreasing in a vicinity of $1/e$.
Thus for every fixed $\nu$ there exists a $\tau_x$ strategy with winning probability strictly higher than $W(\tau_{1/e},\nu)$.
\end{proof}

\begin{assertion} Suppose strategy $\tau$ has cutoffs satisfying $a_k\to x$ and $a_1>0$. Then 
$W(\tau,n)\to -x\log x$ as $n\to\infty.$
In particular, for $x=1/e$ the strategy is asymptotically optimal.
\end{assertion}
\begin{proof} Suppose $0<x<1$ and fix $\delta>0$. 
all cutoffs bigger than some $k$ are within $(x-\delta, x+\delta)$.  Since the number of arrivals on $[0,a_1$ goes to infinity with $n$, the probability to stop before $x-\delta$ approaches $0$. 
The probability to stop within $(x-\delta, x+\delta)$ has asymptotic upper bound $3\delta/x$, since the process of record times converges to Poisson process with intensity  $1/t$. 
Finally, on the event  $\tau>x+\delta$ the strategy coincides with $\tau_{x+\delta}$, and the result follows by sending $\delta\to 0$.
The  extreme values $x=0$ and $1$ are treated similarly.
\end{proof}

Analysis of stopping strategies with  multiple cutoffs is much more involved, but a key idea is seen from the classic best choice problem.
Fix $k$. For $n<k$ the value of cutoff $a_k$ for $n$-optimality does not matter. As $n$ increases, it is optimal to keep $a_k=0$ as long as $k\geq d^*(n)$, switching to $a_k=1$ for $n$ about $ke$.
The next lemma quantifies the winning probability as the  cutoffs vary between the extreme positions.

Let $w_n(\boldsymbol{a})$ be $W(\tau, n)$ for $\tau$ the strategy with cutoffs $\boldsymbol{a}=(a_1,a_2,\cdots)$. Note that $w_n$ only depends on $a_1,\cdots,a_n$.

\begin{lemma} Suppose $a_k\leq \min(a_1,\cdots,a_{k-1})$, then for $n\geq k$
\begin{equation}\label{diff-f}
\frac{\partial w_n(\boldsymbol{a})}{\partial a_k}={n-1\choose k-1}  a_k^{k-1} (1-a_k)^{n-k}\,[h(k+1,n)-1].
\end{equation}
\end{lemma}

\begin{proof} Compare strategy $\tau$ with cutoffs $\boldsymbol{a}$ and $\tau'$ with $k$th cutoff changed to $a_k-\delta$, and all other being unaltered.
The strategies  make different choices only if there is a record with index $k$ and arrival time  
  between $a_k-\delta$ and $a_k$.
In that event both strategies do not stop before $a_k-\delta$, because the preceding records
have  index below $k$ and 
 by the assumption on cutoffs could not be accepted. Then, $\tau'$ stops and wins with probability $k/n$, while $\tau$ picks the next record 
(check the acceptance condition!) and wins with probability $s(k+1,n)$. 
The formula follows from the formula for the density of the $k$th uniform order statistic and (\ref{hh}).
\end{proof}

Note that the only zero derivative in  (\ref{diff-f})
 is $\partial w_2/\partial a_1=0$, with all other being sign-definite, a property which will turn important for constructing constrained dominating strategies.
For simplicity of exposition, we choose all cutoffs to the left of $1/e$.

\begin{assertion}\label{Prop6}
For every finite set  of integers  $Z$  there exists a strategy strongly dominating  the $1/e$-strategy on $Z$. 
\end{assertion}
\begin{proof}
Fix any large enough $n_1$, and set iteratively $n_j=d^*(n_{j-1})-1$ until $n_m=0$ for some $m$, resolving the ambiguity for $n=2$ arbitrarily. 
We introduce a sequence of strategies $\tau_1,\cdots, \tau_m$, where $\tau_j$ dominates the $1/e$-strategy for $n\leq n_1$, while dominating strongly for $n_{j+1}<n\leq n_j$.

Define $\tau_1$ by
setting $a_n=1$ for $n_2<n\leq n_1$, and leaving the remaining cutoffs at $1/e$. 
According to (\ref{diff-f}) this gives some improvement.
Over this range of $n$, let $\alpha_1$ be the minimum advantage over the $1/e$-strategy, $W(\tau_1, n)-W(\tau_{1/e},n)$. For smaller values of $n$ the difference is zero. 

Inductively, with $\tau_{j-1}$ defined, let $\alpha_{j-1}$ be the minimum advantage of this strategy for $n>n_{j}$.
To define $\tau_{j}$, set cutoffs $a_n$ for  $n_{j+1}<n\leq n_j$ equal to some value $a_{n_j}$, leaving other cutoffs same as for $\tau_{j-1}$.
As $a_{n_j}$ decreases from $1/e$, a strict advantage for  $n_{j+1}<n\leq n_j$ is gained, but disadvantage for $n>n_j$ increases. Choose 
$a_{n_j}$ in such a way that $a_{n_{j-1}}< a_{n_j}<1/e$ and the advantage for $n>n_j$ is still at least $\alpha_{j-1}/2$, say.
\end{proof}

If the sequence of cutoffs is nonincreasing,  there is an integral counterpart of (\ref{diff-f})

\begin{equation}\label{int-f}
w_n({\boldsymbol a})=\sum_{k=0}^{n-1}  \pi_{n,k}({\boldsymbol a}) \, s(k+1,n),
\end{equation}
where for $0<k<n$
\begin{equation}
\pi_{n,k} ({\boldsymbol a}) =n{n-1\choose k-1} \int_{a_{k+1}}^{a_k}   x^{k-1}(1-x)^{n-k}{\rm d}x+   
  {n\choose k} a_{k+1}^k(1-a_{k+1})^{n-k},
\end{equation}
and $\pi_{n,0}({\boldsymbol a})=(1-a_1)^n$.

More  detailed discussion of the cutoff monotonicity and derivation of (\ref{int-f}) will appear elsewhere.
In general, however,  a $\nu$-optimal strategy need not have nonincreasing cutoffs, as the next example demonstrates. The possibility of 
such irregularity (analogous to stopping islands in \cite{Presman}) was mentioned in \cite{BrussSam1}, p. 827.

\begin{example}{\rm  Consider a two-point distribution $\nu_{10}=\nu_{100}=1/2$. Since $d^*(10)=4$ and $d^*(100)=38$, a $\nu$-optimal strategy will have 
$1=a_1=a_2=a_3>
a_{4}>\cdots>a_{10}>0$. Then 
$a_{11}=\cdots=a_{37}=1, a_{38}=\cdots,a_{100}=0$,
because after the $11$th arrival is observed, the only remaining option is $N_1=100$.
}
\end{example}

\section{Alternatives to the $1/e$ strategy}

\subsection{Dominance for $n>1$}
\label{n>1}
For $0<x< a_1\leq1$, 
define $\tau$ by choosing cutoffs
$a_1=1$ and $a_k=x$ for $k>1$. In the event $N_{x}>0$ the strategy coincides  with $\tau_{x}$. If $N_{x}=0$ the strategy 
selects the second record (if any).

Suppose $N_{x}=0$. Then $\tau$ skips the first arrival and wins with probability $s(2,n)$, while $\tau_{x}$ stops and wins with probability $s(1,n)$.  
From (\ref{int-f})

$$
{W(\tau,n)-W(\tau_{x},n)}=[ (1-x)^n-(1-a_1)^n](s(2,n)-s(1,n)),
$$
where
$$
s(2,n)-s(1,n)=
\begin{cases}
-1, ~~~n=1,\\
0,~~~n=2,\\
h(3,n)/n, ~~~n>2.
\end{cases}
$$
Thus $\tau$ dominates $\tau_x$ for $n>1$ (that is, on the set $Z=\{2,3,\cdots\}$), and strongly dominates for $n>2$. For $n$ large, the advantage is 
asymptotic to $(1-x)^{n}(\log n)/n$.

In particular, setting $a_1=1, x=1/e$  the strategies  compare
for $n=1$ as $0$ to $1-1/e$, for $n=2$ both win with probability $\frac{1}{2}(1-1/e^{2})$,
and for $n=3$ the advantage over the $1/e$-strategy is $\frac{1}{6}(1-1/e)^2=0.066\cdots$.

Setting $a_1=1-1/e, x=1/e$  yields a strategy that dominates $\tau_{1/e}$ for $n\geq 2$, dominates strongly for $n>2$ and for $n=1,2$ has the same winning 
probability as the $1/e$-strategy. 
Clearly, this strategy is minimax.

We see that if Pierre is prohibited to play $n=1$, or if the game starts with an item that sets a standard but is not choosable (like in the models in \cite{Gaver, Campbell}) then there is a good reason to discard the $1/e$-strategy.

\subsection{Dominance for $n>2$}
We introduce next non-Markovian strategies to strongly dominate some $x$-strategies (including the $1/e$-strategy) for $n>2$.
 A key observation is that 
$p_n$'s for $n>2$ are strictly increasing on $[0,c]$, where $c=2-\sqrt{3}=0.267\cdots$ is the maximum point of $p_3$.

Fix $x>c$ and choose  $y$ in the range 
$$\left(\frac{x-c}{1-c}\right)_+<y\leq x.$$  
For $x'=c(1-y)+y$ let

$$
\tau= \begin{cases}  \tau_x, ~~{\rm if~~} N_y>0,  \\
\tau_{x'}, ~~{\rm if~~} N_y=0.
\end{cases}
$$
Given $N_y=0$, the arrivals sequence can be identified with order statistics from the uniform distribution on $[y,1]$, hence $\tau$ wins with 
probability $p_n({c})$, and $\tau_x$ with probability $p_n((x-y)/(1-y)$, which gives
$$W(\tau, n)-W(\tau_x,n)=(1-y)^n \left[p_n(c)-p_n\left(\frac{x-y}{1-y}\right)\right],$$
which is strictly positive for $n>2$.

In particular, for $x=1/e$ choosing $y=x$ the advantage over the $1/e$ strategy for $n>2$ becomes
$$(1-1/e)^n (|c\log c|-1/n),$$
where $|c\log c|=0.352\cdots$. The surplus results  from the event that the first arrival appears after $1/e$. The exponential factor can be increased by tuning parameter $y$, for instance 
taking $y=1/(2e)$ gives 
$$W(\tau, n)-W(\tau_{1/e},n)=\left(1-\frac{1}{2e}\right)^n \left[p_n(c)-p_n\left(\frac{1}{2e-1}\right)\right],$$
where for $n\to\infty$
$$ p_n(c)-p_n\left(\frac{1}{2e-1}\right)      \to 0.017\cdots.$$
So the advantage is of  higher order than for the strategies in Section \ref{n>1}.

\subsection{Minimax strategies}

A class of minimax strategies is obtained by shifting the cutoffs to the right, the direction opposite to that used
in Proposition \ref{Prop6}.

Formula (\ref{diff-f}) says that it is advantageous to shift cutoff  $a_k$ to the right, provided  the number of items is big enough to satisfy $d^*(n)>k$;
that is, for $n$ such that the $n$-optimal strategy would stop at record with
index $k$. 
Since $\tau_{1/e}$ has winning probability well above $1/e$ for small $n$,
there is some room to trade the winning chance for smaller $n$ against some advantage for larger.

In Section \ref{n>1} we deformed the $1/e$-strategy by chosing $a_1=1/e$, thus  reducing $W(\tau_{1/e},1)= 1-1/e=0.632\cdots$  to $1/e=0.367\cdots$, but gaining for $n\geq 2$ (strictly
for $n>2$).

More generally, think of $1/e$ as the initial position of all cutoffs.
Choose $1/e< a_1\leq 1-1/e$, then 
 increase $a_2$ to a position so as to keep $a_2<a_1=1/e$ and the winning probability for $n=2,3,4$ strictly above $1/e$.
Since $d^*(5)=3$ this improves
 the chances for $n>5$. Then $a_3$ can be increased  subject to similar constraints for $n=3,4,5,6,7$ ($d^*(8)=4$), 
and so on. Every step gives a minimax strategy that dominates $\tau_{1/e}$ strongly on $\{n, n+1,\cdots\}$ but is dominated by $\tau_{1/e}$ on the first $n$ integers.
Infinitely many steps result in a minimax strategy dominated by $\tau_{1/e}$.

Computation with (\ref{int-f}) shows that 
for the first three cutoffs it turns that the only active constraints is the winning probability for $n=1,2,3$. Pushing these to the limits, we get  
\begin{eqnarray*}
a_1&=&1-\frac{1}{e}=0.632\cdots,\\
a_2&=&\left(1-\frac{2}{e}  \right)^{1/2}=0.514\cdots,\\
a_3&=&\frac{1}{2^3}\left[2\left(1-\frac{2}{e} \right)^{3/2}-\frac{1}{e^3}   \right]^{1/3}=0.480\cdots,
\end{eqnarray*}
which defines a strategy which has
$$w_n(a_1,a_2,a_3, 1/e, 1/e,\cdots)=1/e=0.367\cdots, ~~~{\rm for ~}n=1,2,3,$$ 
exactly, and beats $\tau_{1/e}$ starting from $n=5$
(though the theory above guaranteed improvement for $n\geq 8$):
\vskip0.2cm
  \begin{tabular}{ | l | c | c|c|c|c| c|c|c|c|c|}
    \hline
   $ n$ &   4&5&6&7&8&9 & 10&15&20&25 \\ 
\hline
   $ w_n$ &  0.373    &0.379     &0.383 & 0.385 &0.384 &0.382 &0.380&0.370&0.368&0.367\\
 \hline
$p_n(1/e)$& 0.376&0.371&0.369&0.368&0.368&0.368&0.368 & 0.367&0.367&0.367\\
\hline

  \end{tabular}

\paragraph{Remark} It is not hard to see that a strategy with nonincreasing cutoffs allocated on both sides
from $1/e$ cannot dominate $\tau_{1/e}$. 
Evaluation of such strategies is  more difficult to be used for verifying the dominance conjecture stated in Introduction.

\section{Counter-examples to monotonicity relative to\\ the stochastic order}

The winning probability $W(\tau,n)$ decreases with $n$ for $\tau=\tau_x$, but may fluctuate in general (the last example).
Also, note that for Markovian strategy with nonincreasing cutoffs, $\tau$ is nonincreasing with $n$, as is seen by coupling different sizes.

The monotonicity of winning probability  holds for the $d^*$-strategy, 
and  this is very intuitive as the strategy depends on $n$ and the choice problem  becomes harder.
For random sample size the natural sense of monotonicity is relative to the (partial) stochastic order $\succ$ on distributions. Recall that $\nu'\succ \nu$ if $\nu'$ has heavier tails.

Theorem 2.3 from \cite{BrussSam1}, when adjusted to the best-choice context, asserts that the relation $\nu'\succ \nu$  implies $\max W(\,\cdot\,,\nu')\leq  \max W(\,\cdot\,,\nu)$, in line with the 
fixed number of items case. We will disprove this by counter-examples.
A minor delicate point is that in \cite{BrussSam1}, p. 825, the payoff of non-stopping in case $n=1$ is  1
(note that in this setting  Pierre will never play $n=1$, and the $1/e$-strategy will be dominated, see Section \ref{n>1}).
The monotonicity claim is re-stated in \cite{BrussSam2} under the opposite convention that the non-stopping in the $n=1$ case is assessed as $0$.
The first, simpler, example to follow is a counter-example under the second convention, and the second works for both.

\begin{example} {\rm Suppose Pierre plays a two-point distribution $\nu_1=1-p, \nu_4=p$.
The distribution is strictly increasing in the stochastic order as $p$ increases from $0$ to $1$.
   
We have $d^*(1)=1, d^*(4)=2$. If the first arrival is not chosen, then regardless of the time
the best way to proceed  it to stop at the next record (with a hope that
   $n=4$).
Thus  $a_2=a_3=a_4=0$, and the strategies failing the condition can be discarded by dominance.
We leave to the reader a rigorous proof that the optimal acceptance region $A_1$ is an interval, but this is 
intuitively obvious in this simple situation.
Thus the only indeterminate of the stopping strategy is the cutoff $a_1$.
Given the first item comes before the cutoff, Alice proceeds with  
$2$-strategy, otherwise with $1$-strategy. 
Changing the variable as $b=1-a_1$ for shorthand, the total winning 
probability is computed as
\begin{eqnarray*}
w(b,p):=(1-p)b+p \left((1-b^4)\,s(2,4) +\frac{b^4 }{4}\right)&=&\\
(1-p)b+\frac{p}{24}(11- 5b^4), ~~~(b,p)\in[0,1]^2,
\end{eqnarray*}]
where we used
$s(2,4)=11/24$. For $p\leq 6/11$ this  increases in $b$, hence the best response has $a_1=1-b=0$ and  the strategy always stops  at the first arrival.
For $0<p<6/11$ the function $w(p,\cdot)$ has a single mode inside $(0,1)$. The saddle point is
$$b^*=0.449\cdots,~~~~  p^*=0.929\cdots,$$
where $b^*$ is a root of
$$-b+\frac{1}{24}(11-5b^4)=0.$$
This is an equaliser, hence
$$w(b^*,p)=b^*, ~~~p\in [0,1].$$
Thus if Pierre plays $\nu_1=1-p^*,\nu_4=p^*$, Alice can only achieve $b^*$, and so for $p\in [p^*,1]$
$$b^*=w(b^*,p)< \max_b w(b,p)\leq w(0,1)= s(2,4)=\frac{11}{24}=0.458\cdots,$$ 
where the second term increases in $p$.
The conclusion is that a lottery on $\{1,4\}$ may turn for Alice less favourable than certain $n=4$.
}
\end{example}

\begin{example} {\rm Suppose Pierre plays $\nu_3=1-p, \nu_6=p$. Since $d^*(3)=2, d^*(6)=3$, Alice will  play $a_1=1, a_3=a_4=a_5=a_6=0$ and some 
$a_2$. Let $b=1-a_2$ for shorthand, $\beta(k,n):={n\choose k}(1-b)^k b^{n-k}$. If $n=3$ 
Alice wins with probability
$$f_1(b):= [\beta(0,3)+\beta(1,3)]s(2,3)+ [(\beta(2,3)+\beta(3,3)]s(3,3),$$
and if $n=6$ with probability 
$$f_2(b):= [\beta(0,6)+\beta(1,6)]s(2,6)+ (\beta(2,6)+ \beta(3,6)+\beta(4,6) +\beta(5,6)+\beta(6,6) ]          s(3,6).$$
The game on the unit square has the payoff `matrix'
$$w(b,p):= (1-p)f_1(b)+p f_2(b).$$ 
With the experience of the previous example, we look for an equalising strategy.
Equation $f_1(b)=f_2(b)$ becomes
$$\frac{17}{72}\,b^6- \frac{17}{60}\, b^5  +  \frac{1}{3}\,b^3   -  \frac{1}{2} \,b^2   + \frac{17}{180}=0,$$
and has a unique suitable root $b^*=0.520\cdots$, with $w(b^*,p)=0.421\cdots$ for all $p\in[0,1]$.

To find a best response to $p$, we observe that $w(b,p)$ is maximal at $b(p)=1$ for $p\leq 0.413\cdots$,
and at 
$$b(p)=\left(\frac{12(1-p)}{17p}\right)^{1/3} ~~~{\rm for~} p>0.413\cdots,$$
as found by solving
$$\frac{\partial w(b,p)}{\partial p}= (1-p)b(1-b)-\frac{17}{12}\,p \,b^4(1-b)=0.$$
Minimising $w(b(p),p)$ returns $p^*=0.833\cdots$, so $(b^*,p^*)$ is a saddle point.
The $\nu$-optimal winning  probability $w(b(p),p)$ is strictly increasing on $[p^*,1]$. But the larger $p$, the larger $\nu$  stochastically.
}
\end{example}

\paragraph{Remark}  A gap in  \cite{BrussSam1} (Theorem 2.3)
appeared in the short argument on p. 827, where dependence of the stopping and continuation risks on $N_t$ was ignored. 
The asserted parallel with  Theorem 2.1 of \cite{BrussSam1} is not relevant here, as the result  concerns  strategies that rely  decisions
 on the arrival times  and the relative ranks only, while even the classic  $d$-strategies (embedded in continuous time) are not of this kind.
Nevertherless, the implication
$\nu'\succ \nu\Rightarrow\max W(\,\cdot\,,\nu')\leq  \max W(\,\cdot\,,\nu)$ 
does hold  under the additional assumption that  $\nu'$ is a convolution of $\nu$ with another distribution on nonnegative integers; in that 
case the proof follows by coupling exactly as in the fixed sample size problems \cite{BG, Gianini-P}.

\end{document}